\documentclass[reqno]{amsart}
\usepackage[utf8x]{inputenc}  

\usepackage{setspace}
\usepackage[left=3.1cm, right=3.1cm, bottom=4cm]{geometry}       

\usepackage{graphicx} 
\usepackage{pgf,tikz}
\usetikzlibrary{arrows}
\usepackage{amssymb}
\usepackage{pdfsync}
\usepackage{mathrsfs}
\usepackage{hyperref} 
\usepackage{verbatim} 
\usepackage{epstopdf}
\usepackage{bbm} 
\usepackage[colorinlistoftodos,prependcaption,textsize=tiny]{todonotes}
\usepackage{xargs}
\usepackage{bm}

\def\today{\ifcase\month\or
  January\or February\or March\or April\or May\or June\or
  July\or August\or September\or October\or November\or December\fi
  \space\number\day, \number\year}
\DeclareMathOperator{\sgn}{\mathrm{sgn}}

 \newtheorem{theorem}{Theorem}

 \newtheorem{corollary}[theorem]{Corollary}
 \theoremstyle{definition}

 \theoremstyle{remark}

 \newcommand{\C}{\mathbb{C}}
 \newcommand{\R}{\mathbb{R}}
  
 \newcommand{\N}{\mathbb{N}}
 
 \newcommand{\Z}{\mathbb{Z}}

 \newcommand{\ds}{\text{\rm d}s}

 \newcommand{\du}{\text{\rm d}u}

 \newcommand{\dx}{\text{\rm d}x}

\begin{document}

\title[Monotone extremal functions]{Monotone extremal functions and \\ the weighted Hilbert's inequality}
\author[Carneiro and Littmann]{Emanuel Carneiro and Friedrich Littmann}
\date{\today}

\address{ICTP - The Abdus Salam International Centre for Theoretical Physics, 
Strada Costiera, 11, I - 34151, Trieste, Italy.}
\email{carneiro@ictp.it}
\address{Department of mathematics, North Dakota State University, Fargo, ND 58105-5075.}
\email{friedrich.littmann@ndsu.edu}

\allowdisplaybreaks
\numberwithin{equation}{section}

\maketitle

\begin{abstract} In this note we find optimal one-sided majorants of exponential type for the signum function subject to certain monotonicity conditions. As an application, we use these special functions to obtain a simple Fourier analysis proof of the (non-sharp) weighted Hilbert-Montgomery-Vaughan inequality.
\end{abstract}

\section{Introduction} An entire function $F: \C \to \C$ is said to be of {\it exponential type} if 
$$\tau(F) := \limsup_{|z| \to \infty} \,|z|^{-1}  \log |F(z)|  < \infty.$$
In this case, the number $\tau(F)$ is called the {\it exponential type} of $F$. An entire function $F:\C \to \C$ is said to be {\it real entire} if it its restriction to $\R$ is real-valued. In this note we solve the following extremal problem with a monotonicity constraint.

\begin{theorem} \label{Thm1}
Let $F:\C \to \C$ be a real entire function such that:
\begin{itemize}
\item[(i)] $F$ has exponential type at most $2\pi$;
\item[(ii)] $F(x) \geq \sgn(x)$ for all $x \in \R$;
\item[(iii)] $F$ is non-decreasing on $(-\infty, 0)$ and non-increasing on $(0,\infty)$.
\end{itemize} 
Then 
\begin{equation}\label{Ext_ineq}
\int_{-\infty}^{\infty} \big\{F(x) - \sgn(x)\big\}\,\dx \geq 2.
\end{equation}
Moreover, there exists a unique real entire function $M:\C \to \C$ verifying properties {\rm (i)}, {\rm (ii)} and {\rm (iii)} for which the equality in \eqref{Ext_ineq} holds. This function is given by
\begin{equation}\label{20230615_16:01}
M(z) = -2  \int_{-\infty}^z  \frac{\sin^2 \pi s}{\pi^2 s\,(s+1)^2}\,\ds - 1.
\end{equation}
\end{theorem}
\noindent {\sc Remark:} The integral in \eqref{20230615_16:01} is understood to be along the path $(-\infty, 0] \cup [0,z]$, where the latter is the line segment connecting $0$ to $z$.

\smallskip

\begin{figure} 
\includegraphics[height=4.5cm]{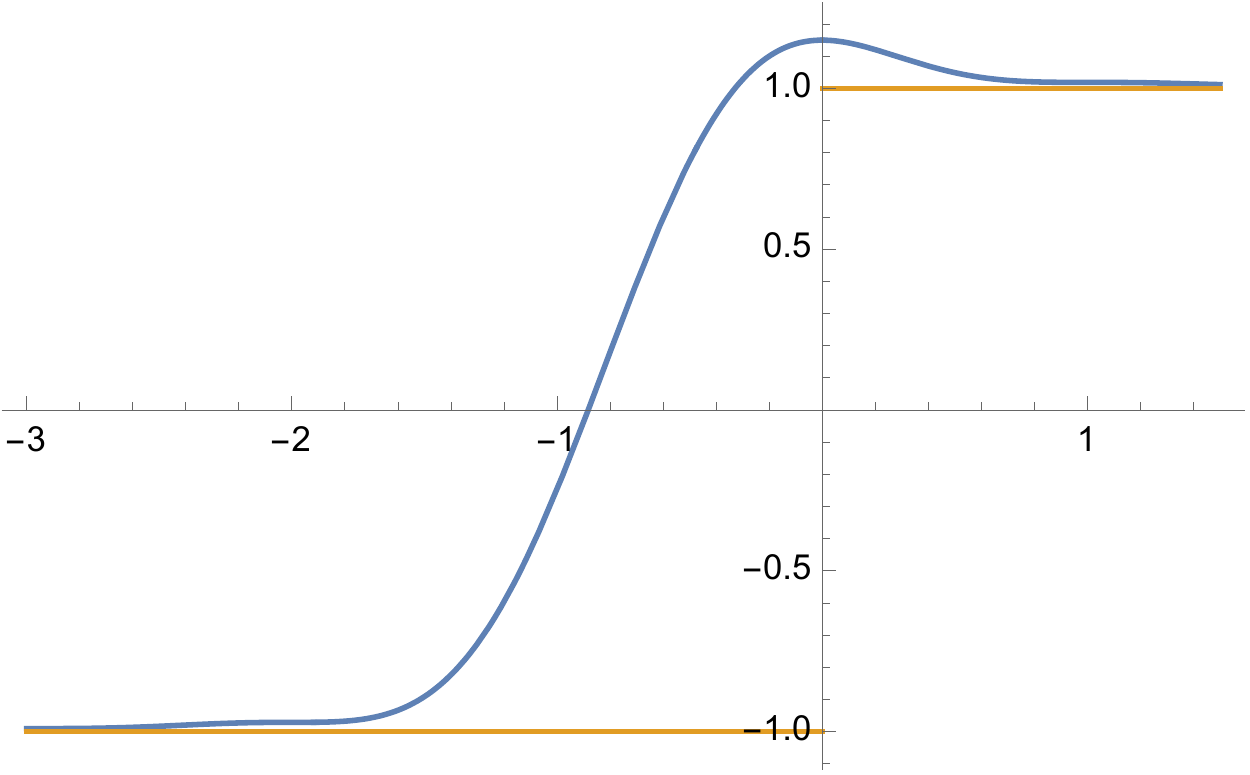} \qquad 
\includegraphics[height=4.5cm]{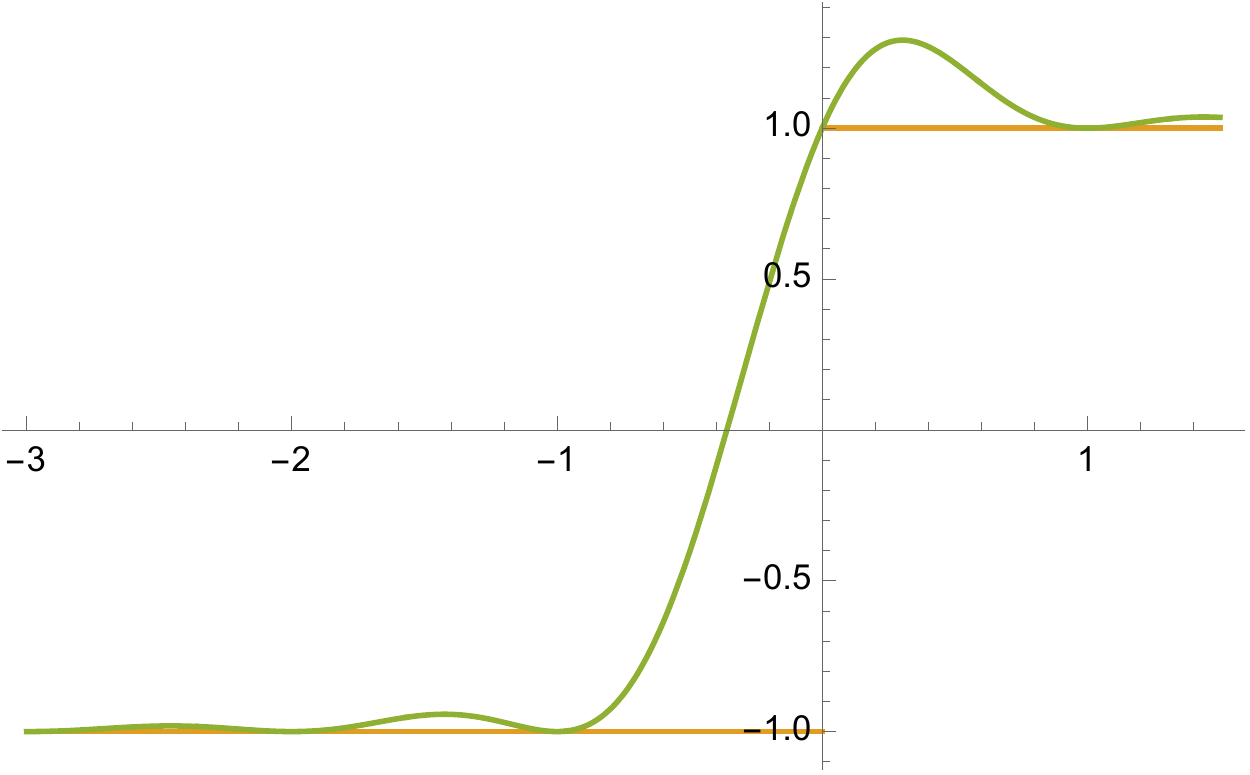} 
\caption{The monotone extremal majorant $M(x)$ on the left, and the classical Beurling majorant $B(x)$ on the right.}
\label{figure1}
\end{figure}

Without the monotonicity constraint (iii) in Theorem \ref{Thm1}, this problem was solved by Beurling in the late 1930's, and the value of the minimal integral on the right-hand side of \eqref{Ext_ineq} is actually equal to $1$\,; see Vaaler's classical survey \cite{V} on the subject. The unique extremal function in this case is
\begin{equation}\label{20230615_20:05}
B(z) = \left(\frac{\sin \pi z}{\pi}\right)^2\left(\sum_{n=0}^{\infty} \frac{1}{(z-n)^2} - \sum_{m=-\infty}^{-1} \frac{1}{(z-m)^2} + \frac{2}{z}\right).
\end{equation}
See Figure \ref{figure1} for the plots of these functions on $\R$.

\medskip

As an application of Theorem \ref{Thm1} we revisit the following result of Montgomery and Vaughan \cite{MV}.
\begin{corollary}[Weighted Hilbert-Montgomery-Vaughan inequality]\label{Cor2}
Let $N\in\N$. Let $\lambda_1,\ldots,\lambda_N$ be a set of distinct real numbers and define $\delta_n:=\min\{|\lambda_n-\lambda_m|:m\neq n\}$. If $a_1,\ldots,a_N\in\C$ then
\begin{align}\label{weighthilbert}
\left|\sum_{\substack{m,n =1 \\ m\neq n}}^N\frac{a_m\overline{a}_n}{ (\lambda_m-\lambda_n)}\right|\le C \, \sum_{n=1}^{N} \frac{|a_n|^2}{\delta_n}
\end{align}
holds with $C = 2\pi$. 
\end{corollary}

Inequality \eqref{weighthilbert} has a long history. In the case $\lambda_m = m$, inequality \eqref{weighthilbert} with constant $C=2\pi$ was first proved by Hilbert. This was later improved by Schur \cite{Schur}, who obtained the sharp constant $C=\pi$ on the right-hand side. The equally-spaced case of \eqref{weighthilbert} (i.e. when the $\{\delta_n\}_{n=1}^N$ on the right-hand side are replaced by a uniform $\delta$) with the sharp constant $C=\pi$ was established by Montgomery and Vaughan in \cite{MV} with a spectral analysis approach, and by Vaaler \cite{V} with a Fourier analysis approach based on Beurling's extremal functions. The general weighted case was first proposed by Montgomery and Vaughan \cite{MV}, who proved inequality \eqref{weighthilbert} with constant $C = \frac{3}{2}\pi$. This was later improved by Preissmann \cite{P} who obtained
$$C = \sqrt{1 + \tfrac{2}{3} \sqrt{\tfrac{6}{5}}} \ \pi = (1.3154...)\,\pi\,,$$ 
currently the best known bound in the literature. Selberg privately reported to Montgomery a proof of \eqref{weighthilbert} with constant $C= 3.2$, but the ideas of such a proof were never made public. It is conjectured that \eqref{weighthilbert} should hold with constant $C=\pi$, and this has been an open problem since 1974. 

\smallskip

Our contribution in this application is to provide, for the first time, a Fourier analysis proof of the weighted Hilbert-Montgomery-Vaughan inequality. Such a proof turns out to be simple, with the caveat of giving a slightly worse constant $C=2\pi$. The previous proofs of Montgomery and Vaughan \cite{MV} and of Preissmann \cite{P} live within the realm of linear algebra, relying on an intricate series of estimates to directly bound the largest eigenvalue of the associated hermitian matrix. 

\medskip

Weighted inequalities like \eqref{weighthilbert} have many applications in number theory, e.g. \cite{M1} and \cite{MV2}. Other works related to the weighted Hilbert-Montgomery-Vaughan inequality include \cite{Li} and \cite{Y}.

\section{A Fourier analysis proof of the weighted \\ Hilbert-Montgomery-Vaughan inequality}

In this section we assume the validity of Theorem \ref{Thm1} and prove Corollary \ref{Cor2}.

\subsection{Proof of Corollary \ref{Cor2}} Let $\psi(x):=M(x)-\text{sgn}(x)$. Throughout this proof we use the notation $\psi_{\delta}(x):=\psi(\delta x)$, for $\delta>0$. By construction, $\psi\in L^1\cap L^2(\R)$, and we denote its Fourier transform on the real line by
\begin{equation*}
\widehat{\psi}(t):=\int_{-\infty}^{\infty} \psi(x) \,e^{-2\pi ixt}\,\dx.
\end{equation*}
We remark that
\begin{align*}
\widehat{\psi_\delta}(t)=\delta^{-1}\widehat{\psi}(\delta^{-1}t),
\end{align*}
and hence, by the Paley-Wiener theorem,
\begin{align*}
\widehat{\psi_\delta}(t)= - (\pi it)^{-1}\text{ for }|t|\ge \delta.
\end{align*}
Reorder the sequence $\{\lambda_n\}_{n=1}^N$ so that $\delta_1\ge \delta_2\ge\ldots\ge\delta_N>0$. Then, evidently, $|\lambda_m-\lambda_n|\ge \delta_{\min(m,n)}$. We adopt the convention $\psi_0\equiv 0$. From the monotonicity condition we note that $\psi_{\delta_{j}}(x) \geq \psi_{\delta_{j-1}}(x)$ for all $x \in \R$ and $j=1,2, \ldots, N$. Hence
\begin{align}
0 & \leq \sum_{j=1}^N \int_{-\infty}^{\infty} \big[\psi_{\delta_{j}}(x)-\psi_{\delta_{j-1}}(x)\big]\left|\sum_{m=j}^N a_m e^{-2\pi i\lambda_m x}\right|^2 \,\dx \label{Main_Chain} \\
& = \sum_{j=1}^N\sum_{m,n=j}^N a_m\overline{a}_n\big[\widehat{\psi_{\delta_{j}}}(\lambda_m-\lambda_n)-\widehat{\psi_{\delta_{j-1}}}(\lambda_m-\lambda_n)\big] \nonumber \\
&= \sum_{m,n=1}^N  a_m\overline{a}_n\sum_{j=1}^{\min(m,n)}\big[\widehat{\psi_{\delta_{j}}}(\lambda_m-\lambda_n)-\widehat{\psi_{\delta_{j-1}}}(\lambda_m-\lambda_n)\big] \nonumber \\
& =\sum_{m,n=1}^N a_m\overline{a}_n\,\widehat{\psi_{\delta}}_{\min(m,n)}(\lambda_m-\lambda_n) \nonumber \\
&= - \sum_{\substack{m,n =1 \\ m\neq n}}^N\frac{a_m\overline{a}_n}{\pi i(\lambda_m-\lambda_n)}+ \widehat{\psi}(0)\sum_{n=1}^N\frac{|a_n|^2}{\delta_n}. \nonumber
\end{align}
It then follows that
\begin{align*}
\sum_{\substack{m,n =1 \\ m\neq n}}^N\frac{a_m\overline{a}_n}{\pi i(\lambda_m-\lambda_n)} \leq  \widehat{\psi}(0)\sum_{n=1}^N\frac{|a_n|^2}{\delta_n}.
\end{align*}
The function $-M(-x)$ is a minorant of $\text{sgn}(x)$ which is non-increasing on $(-\infty,0)$, and non-decreasing on $(0,\infty)$. Repeating the above argument with $\varphi(x):=\text{sgn}(x)+M(-x)\geq 0$ yields
\begin{equation*}
- \widehat{\psi}(0)\sum_{n=1}^N\frac{|a_n|^2}{\delta_n} \leq 
\sum_{\substack{m,n =1 \\ m\neq n}}^N\frac{a_m\overline{a}_n}{\pi i(\lambda_m-\lambda_n)},
\end{equation*}
and this concludes the proof of Corollary \ref{Cor2} since $\widehat{\psi}(0) = 2$.

\smallskip

\noindent {\sc Remark}: One might wonder whether these techniques can be used to prove the sharp weighted Hilbert-Montgomery-Vaughan inequality. If one replaces the function $M$ from Theorem \ref{Thm1} by the original Beurling majorant $B$ described in \eqref{20230615_20:05}, and defines $\psi(x):=B(x)-\text{sgn}(x)$ instead, one would need to verify the non-negativity of the corresponding expression appearing in \eqref{Main_Chain}.

\section{Monotone extremal functions: Proof of Theorem \ref{Thm1}}

To simplify the calculations let us consider the analogous extremal problem replacing $\sgn(x)$ by the upper semi-continuous Heaviside function $x_+^0$ (i.e. $x_+^0=1$ for $x\ge 0$, and $x_+^0=0$ for $x<0$). 

\smallskip

We first find necessary conditions that the optimal function must satisfy, and then construct a function that satisfies these conditions. Let $G$ be an entire majorant of $x_+^0$ of exponential type at most $2\pi$ that is non-decreasing on $(-\infty,0)$ and non-increasing on $(0,\infty)$, and such that $G - x_+^0 \in L^1(\R)$. Then $g = G'$ is of one sign on each of the half-lines. Since $\lim_{x\to-\infty} G(x) = 0$, 
it follows that 
\begin{align*}
G(x) = \int_{-\infty}^x g(x) \,\dx.
\end{align*}
Also, from the fact that $\lim_{x\to \infty} G(x) = 1$, it follows that $g \in L^1(\R)$ and 
\begin{align}\label{g_int_1}
\int_{-\infty}^{\infty} g(x) \,\dx = 1.
\end{align}
Since $G - x_+^0 \in L^1(\R)$, Fubini's theorem and \eqref{g_int_1} imply that the following two integrals are finite:
\begin{align}\label{Gneg}
\int_{-\infty}^0 G(x)\,\dx = \int_{-\infty}^0 \int_{-\infty}^x g(u) \,\du\,\dx = -\int_{-\infty}^0 u \,g(u) \,\du,
\end{align}
and
\begin{align}\label{Gpos}
\int_0^\infty \big\{G(x) - 1\big\} \,\dx = \int_0^\infty\left( -\int_{x}^\infty g(u) \,\du\right) \dx = -\int_0^\infty u\, g(u) \,\du.
\end{align}

\smallskip

Define $H:\C \to \C$ by $H(u) = -u \,g(u)$. Then $H$ is a real entire function of exponential type at most $2\pi$ that is non-negative on $\R$. From \eqref{Gneg} and \eqref{Gpos} we find also that $H \in L^1(\R)$. Moreover, since $g$ is entire, $H$ has a zero at the origin of even order at least $2$. It follows by Krein's decomposition\footnote{If $f:\C \to \C$ is a real entire function of exponential type at most $2\pi$, that is non-negative and integrable on $\R$, then there exists $g:\C \to \C$ entire of exponential type at most $\pi$ such that $f(z) = g(z)\overline{g(\overline{z})}$ for all $z \in \C$.} \cite[p. 154]{A} that there exists $h:\C\to \C$ entire of exponential type at most $\pi$ such that
\begin{equation*}
H(z) = z^2 h(z) \overline{h(\overline{z})}
\end{equation*}
for all $z \in \C$. From \eqref{Gneg}, \eqref{Gpos} and Poisson summation (that holds pointwise everywhere since $H' \in L^1(\R)$ by a classical result of Plancherel and P\'{o}lya \cite{PP}, and hence $H$ has bounded variation on $\R$) we have
\begin{align*}
\int_{-\infty}^\infty \big\{G(x) - x_+^0\big\} \,\dx = \int_{-\infty}^\infty H(x) \,\dx = \sum_{n\in\Z} H(n) = \sum_{n\in\Z} |n \,h(n)|^2.
\end{align*}
Another application of Poisson summation, together with \eqref{g_int_1}, yields 
\begin{align}\label{crucial_cond}
-\sum_{n\in\Z} n\,|h(n)|^2 = \sum_{n\in\Z} g(n) = \widehat{g}(0) = \int_{-\infty}^\infty g(x)\,\dx = 1.
\end{align}
Hence
\begin{align}\label{final_inequality}
\int_{-\infty}^\infty \big\{G(x) - x_+^0\big\} \,\dx = \sum_{n\in\Z} |n\,h(n)|^2 \ge \sum_{n\in\Z} |n|\, |h(n)|^2\ge -\sum_{n\in\Z} n\,|h(n)|^2 = 1,
\end{align}
which establishes the desired inequality.

\medskip

In order to have equality in \eqref{final_inequality} we must have $h(n) = 0$ if $n \neq -1,0$. From \eqref{crucial_cond}, this implies that $|h(-1)| =1$. Since $z \,h(z) \in L^2(\R)$ and has exponential type at most $\pi$, the classical Shannon-Whittaker interpolation formula yields
\begin{equation*}
z \,h(z) = h(-1) \frac{\sin \pi z}{\pi (z+1)},
\end{equation*}
which implies that 
\begin{equation*}
g(z) = -z \,h(z) \overline{h(\overline{z})} = -|h(-1)|^2 \frac{\sin^2 \pi z}{\pi^2 z\,(z+1)^2} = - \frac{\sin^2 \pi z}{\pi^2 z\,(z+1)^2}.
\end{equation*}
One can check directly that this $g$ satisfies \eqref{g_int_1} (e.g. via Poisson summation) and that 
\begin{equation}\label{Def_Ext_G}
G(x) = - \int_{-\infty}^x  \frac{\sin^2 \pi u}{\pi^2 u\,(u+1)^2}\,\du
\end{equation}
is indeed a majorant of $x_0^+$ with 
\begin{align*}
\int_{-\infty}^\infty \big\{G(x) - x_+^0\big\} \,\dx = \int_{-\infty}^\infty |u \,g(u)|\,\du  = 1.
\end{align*}
Finally, observe that $G:\R \to \R$ defined by \eqref{Def_Ext_G} is the restriction to $\R$ of the entire function 
\begin{equation}\label{Ent_G}
G(z) = G(0) - \int_0^z \frac{\sin^2 \pi s}{\pi^2 s\,(s+1)^2}\,\ds.
\end{equation}
The integration is over the line segment connecting $0$ to $z$, and the value $G(0)$ is a linear combination of known constants with decimal expansion $G(0) = 1.0749...$ . If $z = x +iy$, it is clear from \eqref{Ent_G} that $|G(z)| \leq C |z|e^{2\pi |y|}$ for some $C>0$, and therefore $G$ has exponential type at most $2\pi$. This concludes the proof in the case of $x_+^0$.

\medskip

\begin{figure}
\includegraphics[height=5cm]{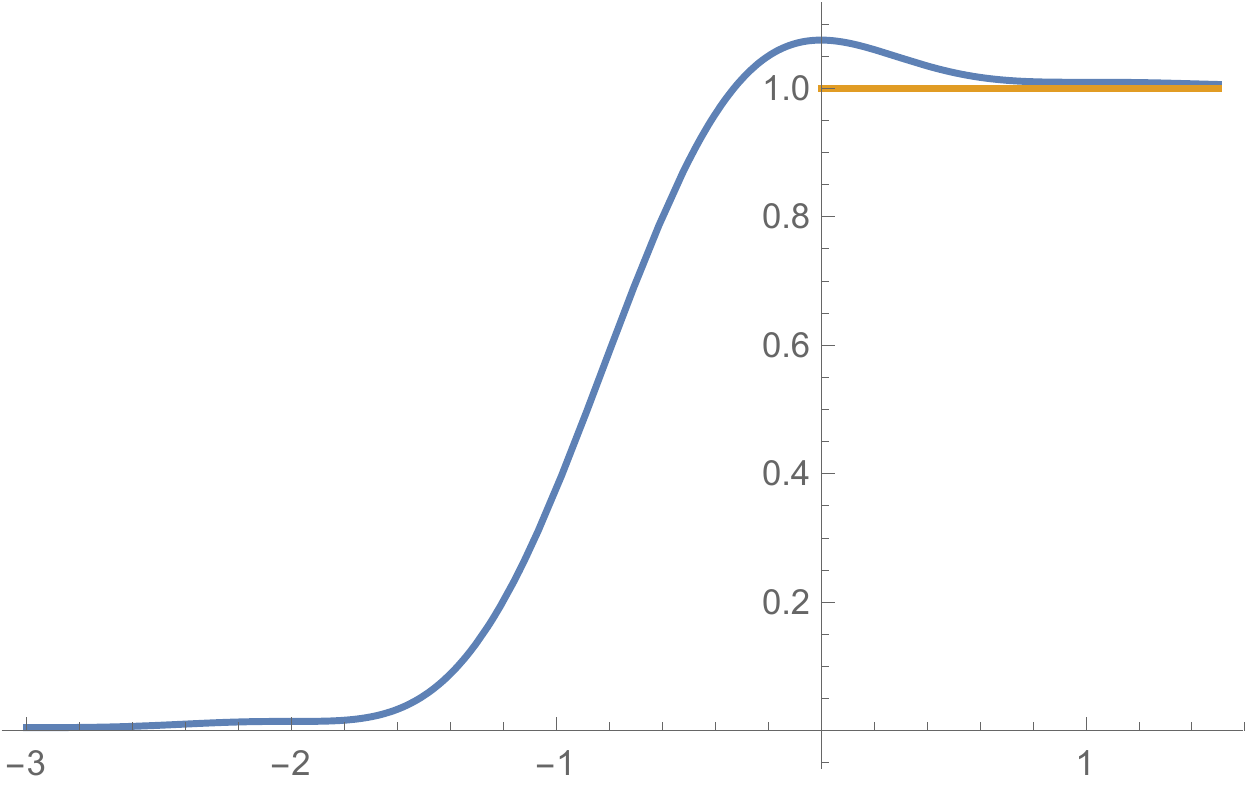}
\caption{The extremal function $G(x)$}
\end{figure}

Naturally, in the case of $\sgn(x)$ our unique extremal function is then $M(z) := 2G(z) -1$.

\section*{Acknowledgements}
We are thankful to Hugh Montgomery and Jeffrey Vaaler for the enlightening discussions on the history of this problem. We are also thankful to Harald Helfgott and Michael Kelly for helpful discussions during the preparation of this note. Finally, we thank the anonymous referees for the helpful remarks in order to improve the presentation.

\end{document}